\documentclass[10pt, leqno]{amsart}
\usepackage{latexsym, amsfonts,mathrsfs, amsmath, amssymb, amscd, epsfig}
\usepackage{mathrsfs}
\usepackage[usenames,dvipsnames]{xcolor}
\usepackage{amsthm}
\usepackage{array}
\usepackage{psfrag}
\usepackage{bm}
\usepackage{setspace}
\usepackage{soul}
\usepackage{graphicx}
\usepackage{pdfpages}
\numberwithin{equation}{section}

\newtheorem{lemma}{Lemma}[section]
\newtheorem{theorem}{Theorem}%[section]
\newtheorem*{question}{Question}
\theoremstyle{definition}

\def\beq#1\eeq{\begin{equation}#1\end{equation}}
\def\balign #1 #2 \ealign{\begin{aligned} #1 #2  \end{aligned} }

\newcommand \alp{\alpha}
\newcommand \eps{\varepsilon}
\newcommand \vphi{\varphi}

\newcommand \Gam{\Gamma}
\newcommand \gam{\gamma}
\newcommand \om{\omega}
\newcommand \tx{\text}
\newcommand \R{\mathbb{R}}
\newcommand \til{\tilde}

\newcommand \der{\partial}
\newcommand \mcl{\mathcal}

\newcommand \ol{\overline}
\newcommand \Om{\Omega}

\newcommand \Gamex{\Gamma_L}
\newcommand \Gamw{\Gamma_w}

\newcommand \rx{\rm{\bf x}}

\newcommand \us{u_s}

\newcommand \ls{l_{\rm s}}

\newcommand \barrhoi{\bar{\rho}_I}

\newcommand \barui{\bar u_{I}}
\newcommand \fsonic{\mathfrak{f}_{\rm sn}}

\newcommand \Tac{\mcl{T}_{\rm acc}}

\newcommand \bJ{\bar{J}}
\newcommand \ubJ{\underbar{J}}

\newcommand \bmxi{\bm{\xi}}

\newcommand \shock{\Gamma_{\rm{shock}}}

\begin{document}
\title[]
{A review on two types of sonic interfaces}

\author{Myoungjean Bae}
\address{Department of Mathematical Sciences, KAIST, 291 Daehak-ro, Yuseong-gu, Daejeon, 43141, Korea}
\email{mjbae@kaist.ac.kr}

\begin{abstract}
In this paper, two examples of sonic interfaces (\cite{bae2013prandtl, bae2019prandtl, bae2023steady, chen2010global, chen2018mathematics}) are presented. The first example shows the case of sonic interfaces as weak discontinuities in self-similar shock configurations of unsteady Euler system.
The second example shows the case of sonic interfaces as regular interfaces in accelerating transonic flows governed by the steady Euler-Poisson system with self-generated electric forces. And, we discuss analytic differences of the two examples, and introduce an open problem on decelerating transonic solution to the steady Euler-Poisson system.
\end{abstract}

\keywords{Keldysh type, transonic, sonic interface, weak discontinuity, regular interface}
\subjclass[2010]{%\AMSMOS
35C06, 35M10, 35M30, 35Q31, 76H05, 76N10}

\date{\today}

\maketitle

%\tableofcontents

\section{A sonic interface as a weak discontinuity}
\label{section:2}

Fix a constant $\eps_0>0$. Given a function $f:[0,\eps_0]\rightarrow \R_+$ with
\begin{equation}
\label{properties-f}
  \|f\|_{C^{1,1}([0,\eps_0])}<\infty,\quad f(0)>0, \quad\tx{and}\quad
  \frac{df}{dx}\ge \om>0\,\forall 0\le x\le \eps_0,
\end{equation}
set
\begin{equation*}
  P_0:=(0,f(0)),
\end{equation*}
and define a domain
\begin{equation}
\label{defi-Q}
  \mcl{Q}_{\eps_0}^f:=\{(x,y): 0<x<\eps_0,\,\,0<y<f(x)\}.
\end{equation}
For each $t\in(0,  f(0))$, set
\begin{equation*}
  \mcl{R}_t:=(0, \frac{\eps_0}{2})\times (0, f(0)-t).
\end{equation*}
Given two constants $a>0$ and $b>0$, and functions $\beta_k\in C(\der \mcl{Q}_{\eps_0}^f\cap \{y=f(x)\})$ for $k=1,2,3$, consider the equation
\begin{equation}\label{2-1}
 (2x-a\psi_x+O_1)\psi_{xx}+O_2\psi_{xy}+
  (b+O_3)\psi_{yy}-(1+O_4)\psi_x+O_5\psi_y=0\quad \mbox{in $\mcl{Q}_{\eps_0}^f$},
\end{equation}
and the boundary conditions
\begin{align}
\label{2-2}
  \psi=0\quad&\mbox{on $\der\mcl{Q}_{\eps_0}^f\cap\{x=0\}$}\\
  \label{2-3}
  \der_y\psi=0\quad&\mbox{on $\der\mcl{Q}_{\eps_0}^f\cap\{y=0\}$}\\
  \label{2-4}
  \beta_1(x,y)\psi_x+\beta_2(x,y)\psi_y+\psi=0\quad&\mbox{on $\der\mcl{Q}_{\eps_0}^f\cap\{y=f(x)\}$}.
\end{align}
In addition, assume that
\begin{equation}
\label{prop-betas}
  \beta_1(x,y)\ge \lambda, |\beta_2(x,y),\, \beta_3(x,y)|\le \frac{1}{\lambda}\quad\tx{on $\der\mcl{Q}_{\eps_0}^f\cap\{y=f(x)\}$}
\end{equation}
for some constant $\lambda>0$.

\begin{theorem}[Theorem 3.1 in \cite{BCF}]\label{theorem:1}
Suppose that a function $\psi:\ol{\mcl{Q}_{\eps_0}^f}\rightarrow \R$ satisfies the following conditions:

\begin{itemize}
\item[(i)] $\displaystyle{\psi\in C^2(\mcl{Q}_{\eps_0}^f)\cap C^{1,1}(\ol{\mcl{Q}_{\eps_0}^f})};$
\item[(ii)] $\psi>0$ in $\mcl{Q}_{\eps_0}^f$;
\item[(iii)] there exist constants $\mu>0$ and $\delta\in(0,1)$ such that
    \begin{equation*}
     -\mu \le \frac{\psi_x(x,y)}{x}\le \frac{2-\delta}{a}\quad\tx{in $\mcl{Q}_{\eps_0}^f$};
    \end{equation*}

\item[(iv)] $\psi$ satisfies \eqref{2-1}--\eqref{2-4}.
\end{itemize}
In addition, assume that the terms $O_i(x,y), i=1,\cdots, 5$ are continuously differentiable, and that there exists a constant $N>0$ satisfying
\begin{equation}
\label{O-terms}
  \begin{aligned}
  \frac{|O_1(x,y)|}{x^2}, \frac{|O_i(x,y)|}{x}\le N\quad&\mbox{for $i=2,\cdots, 5$},\\
  \frac{|DO_1(x,y)|}{x}, |DO_k(x,y)|\le N\quad&\mbox{for $i=2,\cdots, 5$}
  \end{aligned}
\end{equation}
in $\{x>0\}$.
\medskip

Then, the following properties hold:
\begin{itemize}
\item[(a)]$\displaystyle{\forall t\in(0,  f(0)),\,\,\psi\in C^{2,\alp}(\ol{\mcl{R}_t}) \,\,\forall \alp\in (0,1)};$

\item[(b)] $\displaystyle{\psi_{xx}(0,y)=\frac 1a,\quad \psi_{xy}(0,y)=\psi_{yy}(0,y)=0\quad\tx{for all} \quad 0\le y<f(0)};$

\item[(c)] there are two sequences $\{(x_m^{(1)},y_m^{(1)})\}$ and $\{(x_m^{(2)},y_m^{(2)})\}$ in $\mcl{Q}_{\eps_0}^f$ such that
    \begin{equation*}
      \begin{aligned}
      &\lim_{m\to \infty}(x_m^{(1)},y_m^{(1)})=\lim_{m\to \infty}(x_m^{(2)},y_m^{(2)})=P_0
\quad\tx{for $j=1,2$},\\
&\lim_{m\to \infty}\psi_{xx}(x_m^{(1)},y_m^{(1)})=\frac 1a,\quad \lim_{m\to \infty}\psi_{xx}(x_m^{(2)},y_m^{(2)})=0.
      \end{aligned}
    \end{equation*}
\end{itemize}
\end{theorem}

Now we demonstrate an application of Theorem \ref{theorem:1}.

\medskip

An irrotational flow of inviscid compressible polytropic gas is governed by the Euler equation for potential flow:
\begin{equation}
\label{2-6}
  \begin{aligned}
  \der_t\rho+\nabla_{\rx}\cdot (\rho \nabla_{\rx}\Phi)=0,\\
  \der_t\Phi+\frac 12|\nabla_{\rx}\Phi|^2+\frac{\rho^{\gam-1}-1}{\gam-1}=B_0
  \end{aligned}
\end{equation}
for an adiabatic exponent $\gam>1$.
The density and the velocity potential of the flow are represented as $\rho$ and $\Phi$, respectively. And, the term $B_0>0$ represents the Bernoulli constant which is determined by the initial data. For $\theta_w\in (0,\frac{\pi}{2})$, define a symmetric wedge $W$ in $\R^2$ by
\begin{equation}
\label{def-W}
  W:=\{\rx=(x_1, x_2)\in \R^2: |x_2|<x_1\tan \theta_w,\,\,x_1>0\}.
\end{equation}
Suppose that $(\rho, \Phi)$ solves \eqref{2-6} in $\R^2\setminus W$, and satisfies the slip boundary condition
\begin{equation}
\label{2-7}
  \nabla_{\rx}\Phi\cdot {\bf n}_w=0\quad\tx{on $\der W$}
\end{equation}
for the exterior unit normal ${\bf n}$ to $\der W$. Then, for any constant $\alp>0$, it can be directly checked that $(\til{\rho}, \til{\Phi})$ given by
\begin{equation*}
 (\til{\rho}, \til{\Phi})(\rx, t):= (\rho, \frac{1}{\alp}\Phi)(\alp\rx, \alp t)
\end{equation*}
satisfies \eqref{2-6} in $\R^2\setminus W$, and  \eqref{2-7} on $\der W$. Owing to the scaling invariance, one may seek for a self-similar solution in the form of
\begin{equation*}
  (\rho, \Phi)(\rx, t)=(\varrho({\bmxi}), t\Psi(\bmxi))\quad\tx{for}\,\,\bmxi=(\xi_1, \xi_2)\,\,\tx{with}\,\,(\xi_1, \xi_2):=\frac 1t(x_1, x_2).
\end{equation*}
With a pseudo-potential function $\vphi$ given by
\begin{equation*}
\vphi(\bmxi):=-\frac 12|\bmxi|^2+\Psi(\bmxi),
\end{equation*}
one can rewrite \eqref{2-6} as
\begin{equation}
\label{2-8}
  \nabla\cdot(\varrho \nabla\vphi)+2\varrho=0
\end{equation}
with
\begin{equation}
\label{2-10}
  \varrho^{\gam-1}=(\gam-1)B_0+1-(\gam-1)\left(\vphi+\frac 12|\nabla\vphi|^2\right).
\end{equation}

\begin{figure}[htp]
\centering
\begin{psfrags}
\psfrag{r}[cc][][0.7][0]{$S_{\infty}$}
\psfrag{c}[cc][][0.7][0]{$S_0$}
\psfrag{t}[cc][][0.7][0]{$\theta_{\rm w}$}
\psfrag{om}[cc][][0.7][0]{$\Om$}
\psfrag{u}[cc][][0.7][0]{$\Om_0$}
\psfrag{g}[cc][][0.7][0]{$\Gam_0$}
\psfrag{s}[cc][][0.7][0]{$\Gam_{\rm shock}$}
\includegraphics[scale=.5]{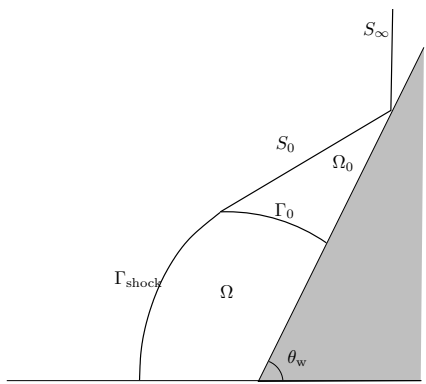}
\caption{Self-similar regular shock reflection ($S_{\infty}$: incoming plane shock, $\ol{S_0\cup\shock}$: reflected shock, $\theta_w\in (\theta_{\rm sonic}, \frac{\pi}{2})$)}
\label{fig1}
\end{psfrags}
\end{figure}

With \eqref{2-8}, the global-in-time existence of weak solutions to \eqref{2-6} of various shock structures are investigated in \cite{bae2013prandtl, bae2019prandtl, chen2010global, chen2018mathematics, chen2008mach, elling2008supersonic}. When a plane shock hits the symmetric wedge $W$ head-on, it is proved in \cite{chen2010global, chen2018mathematics} that if $\theta_w\in(\theta_d, \frac{\pi}{2})$ for the critical angle(=the detachment angle uniquely determined by a shock polar), then \eqref{2-6} has a weak solution of a self-similar regular shock reflection configuration (Fig.\ref{fig1}).

\begin{figure}[htp]
\centering
\begin{psfrags}
\psfrag{c}[cc][][0.7][0]{$S_0$}
\psfrag{t}[cc][][0.7][0]{$\theta_{\rm w}$}
\psfrag{om}[cc][][0.7][0]{$\Om$}
\psfrag{O}[cc][][0.7][0]{$\Om$}
\psfrag{u}[cc][][0.7][0]{$\Om_0$}
\psfrag{g}[cc][][0.7][0]{$\Gam_0$}
\psfrag{h}[cc][][0.7][0]{$\Gam_1$}
\psfrag{s}[cc][][0.7][0]{$\Gam_{\rm shock}$}
\psfrag{d}[cc][][0.7][0]{$S_1$}
\psfrag{v}[cc][][0.7][0]{$\Om_1$}
\includegraphics[scale=.5]{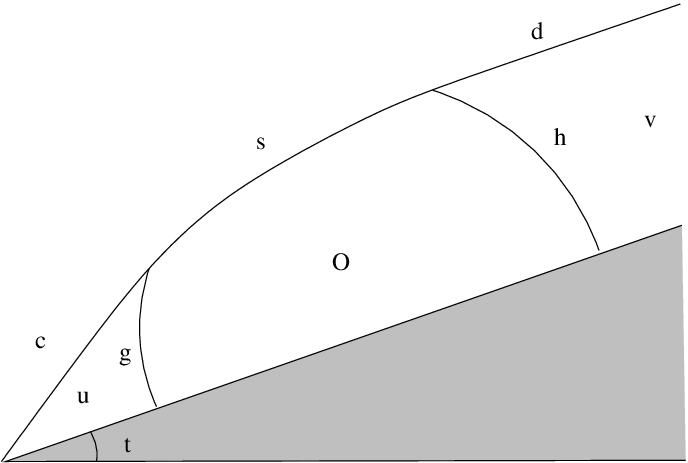}
\includegraphics[scale=.5]{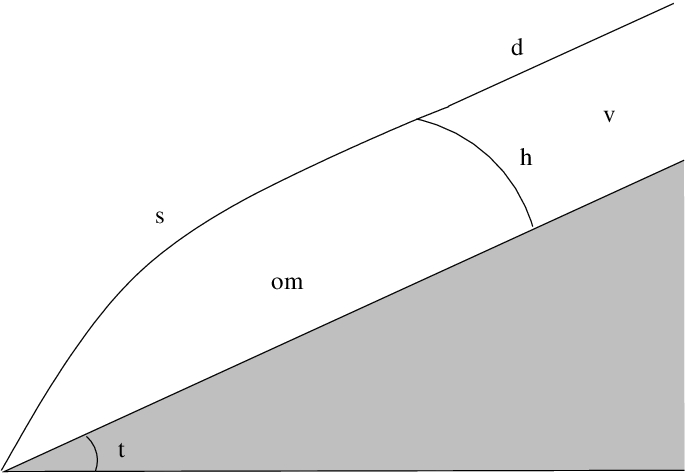}
\caption{Self-similar weak shock configuration past a wedge (Left: $\theta_w<\theta_{{\rm sonic}}$, Right: $\theta_{{\rm sonic}}\le \theta_w<\theta_d$)}
\label{fig2}
\end{psfrags}
\end{figure}

\begin{figure}
\begin{psfrags}
\psfrag{e}[cc][][0.7][0]{$u_2$}
\psfrag{d}[cc][][0.7][0]{$u_1$}
\psfrag{b}[cc][][0.7][0]{$\theta_d$}
\psfrag{a}[cc][][0.7][0]{$\theta_{\rm sonic}$}
\psfrag{c}[cc][][0.7][0]{$\theta_w$}
\psfrag{p}[cc][][0.7][0]{${\bf u}_0$}
\includegraphics[scale=.7]{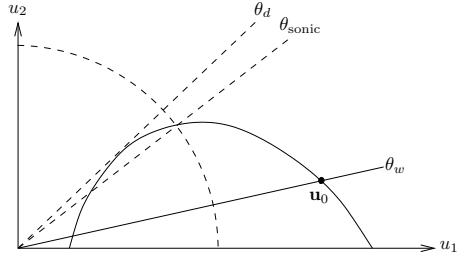}
\caption{The shock polar for steady potential flow equation(${\bf u}_0$: the weak shock state in $\Om_0$ for $0<\theta_w<\theta_{\rm sonic}$)}
\label{fig3}
\end{psfrags}
\end{figure}

When a supersonic flow with a constant density moves horizontally toward the symmetric wedge $W$ at a constant speed, it is proved in \cite{bae2019prandtl, elling2008supersonic} that there exists a global-in-time weak solution of a self-similar weak shock configuration for $\theta_w\in(0, \theta_d)$ (Fig.\ref{fig2}).
In Fig.\ref{fig1} and Fig.\ref{fig2}, the straight shocks $S_0$ and the corresponding downstream state in $\Om_0$ are given as the weak shock state on the shock polar (Fig.\ref{fig3}).
Considering that the incoming supersonic flow has a constant density and a constant velocity, the straight normal shock $S_1$ and the corresponding downstream state in $\Om_1$ (see Fig.\ref{fig2}) can be easily computed from the Rankine-Hugoniot jump conditions.

Since the density $\rho$ in $\Om_0$(Fig.\ref{fig2}) is a constant, it follows from \eqref{2-8} and \eqref{2-10} that the pseudo-potential function $\vphi$ satisfies
\begin{equation*}
\begin{cases}
  \Delta\vphi+2=0\\
  \vphi+\frac 12|\nabla\vphi|^2=\tx{a constant}
\end{cases}\quad\tx{in $\Om_0$},
\end{equation*}
and this yields
\begin{equation}
\label{pseudo potential0}
  \vphi(\bmxi)=-\frac 12|\bmxi|^2+{\bf u}_0\cdot {\bmxi}+k\quad\tx{in $\Om_0$}
\end{equation}
for the constant vector ${\bf u}_0$ given on the shock polar (Fig.\ref{fig3}), and for a constant $k$.

Take the representation of $\varrho$ in terms of $(\vphi, \nabla \vphi)$, directly given from \eqref{2-10}. Substituting the representation $\varrho(\vphi, \nabla \vphi)$ into \eqref{2-8}, we are given with a second order quasi-linear equation
\begin{equation}
\label{2-11}
   \nabla\cdot\left(\varrho(\vphi, \nabla\vphi) \nabla\vphi\right)+2\varrho(\vphi, \nabla\vphi) =0.
\end{equation}
If $\varrho^{\gam-1}-|\nabla\vphi|^2<0$, the equation \eqref{2-11} is hyperbolic, and if $\varrho^{\gam-1}-|\nabla\vphi|^2>0$, it is elliptic.

Back to the uniform state in $\Om_0$(Fig.\ref{fig1}, Fig.\ref{fig2}), let $\rho_0>0$ represent the constant density in $\Om_0$. Then, for the pseudo-potential function $\vphi$ given by \eqref{pseudo potential0}, we have
\begin{equation*}
  \mcl{C}_0:=\{\bmxi: |\nabla\vphi|^2=\varrho_0^{\gam-1}\}=\{|\bmxi-{\bf u}_0|^2=\varrho_0^{\gam-1}\}.
\end{equation*}
The boundary portion $\Gam_0:=\der{\Om_0}\cap \mcl{C}_0$ is called a {\emph{pseudo-sonic arc}} associated with the state of the density $\rho_0$, and the velocity ${\bf u}_0$. In $\Om_0$, the equation \eqref{2-11} is hyperbolic, and the hyperbolicity is degenerate on the pseudo-sonic arc $\der{\Om_0}\cap \mcl{C}_0$. Note that the pseudo-supersonic region $\Om_0$ bounded by a straight oblique shock $S_0$ is shown for $\theta_w\in(\theta_{\rm sonic}, \frac{\pi}{2}]$ in Fig.\ref{fig1}, and for $\theta_w\in(0, \theta_{\rm sonic})$ in Fig.\ref{fig2}. And, this region shrinks to a point as $\theta_w$ approaches to the sonic angle $\theta_{\rm sonic}$, which is uniquely determined by the shock polar associated with the incoming flow state.

The region in which the equation \eqref{2-11} is elliptic($\varrho^{\gam-1}-|\nabla_{\bmxi}\vphi|^2>0$) is indicated as $\Om$. For the admissible solution $\vphi$ constructed in \cite{bae2019prandtl, chen2010global, chen2018mathematics}, one of the essential properties is that there exists a constant $\mu>0$ satisfying
\begin{equation*}
  \frac{|\nabla\vphi(\bmxi)|}{\varrho^{\frac{\gam-1}{2}}(\nabla\vphi(\bmxi),\vphi(\bmxi))}\le 1-\mu\, {\rm dist}(\bmxi, \mcl{C}_0)\quad\tx{for $\bmxi\in \Om$ near $\Gam_0$}
\end{equation*}
as long as the pseudo-supersonic region $\Om_0$ appears (Fig.\ref{fig1}, Fig.\ref{fig2}(left)). 
\medskip

Let $\vphi_0$ be the pseudo-potential function given by \eqref{pseudo potential0}. Define a polar coordinate system by
\begin{equation*}
  \bmxi-{\bf u}_0:=r(\cos\theta, \sin \theta).
\end{equation*}
And, define a new coordinate system by
\begin{equation*}
  (x,y):=\begin{cases}
  (\varrho_0^{\frac{\gam-1}{2}}-r, \theta-\theta_w)\quad&\mbox{for $\bmxi$ near $\Gam_0$ in Fig. \ref{fig1}}\\
  (\varrho_0^{\frac{\gam-1}{2}}-r, \pi+\theta_w-\theta)\quad&\mbox{for $\bmxi$ near $\Gam_0$ in Fig. \ref{fig2}}
  \end{cases}
\end{equation*}
Finally, define
\begin{equation}
\label{def-psi}
  \psi(x,y):=\vphi(\bmxi)-\vphi_0(\bmxi)\quad\tx{in $\Om$ near $\Gam_0$}.
\end{equation}
The admissible solutions constructed in \cite{bae2019prandtl, chen2010global, chen2018mathematics} have the following properties: if $\theta_w\in(\theta_{\rm sonic}, \frac{\pi}{2}]$ in Fig. \ref{fig1}, or if $\theta_w\in(0, \theta_{\rm sonic})$ in Fig. \ref{fig2}, then
\begin{itemize}
\item[(i)] $\exists$ a (small) constant $\eps_0>0$ and a function $f:[0,\eps_0]\rightarrow \R_+$ such that
    \begin{equation*}
      \Om\cap\{{\rm dist}(\bmxi, \Gam_0)<\eps_0\}=\mcl{Q}_{\eps_0}^f
    \end{equation*}
    for the domain $\mcl{Q}_{\eps_0}^f$ defined by \eqref{defi-Q};
\item[(ii)] such a function $f$, representing the curved pseudo-transonic shock $\shock$ (Fig. \ref{fig1}, Fig. \ref{fig2}), satisfies all the properties stated in \eqref{properties-f};
\item[(iii)] the equation \eqref{2-11} is rewritten as
\begin{equation*}
 (2x-(\gam+1)+O_1)\psi_{xx}+O_2\psi_{xy}
 +\left(\frac{1}{\varrho_0^{\frac{\gam-1}{2}}}+O_3\right)\psi_{yy}-(1+O_4)\psi_x+O_5\psi_y=0
\end{equation*}
for the terms $O_j(x,\psi,\psi_x, \psi_y)$ satisfying all the properties stated in \eqref{O-terms};

\item[(iv)] $\psi$ satisfies the boundary conditions \eqref{2-2}--\eqref{2-4} with \eqref{prop-betas} holding;
\item[(v)] $\displaystyle{\psi\in C^2(\mcl{Q}_{\eps_0}^f\cap C^{1,1}(\ol{\mcl{Q}_{\eps_0}^f}))}$;
\item[(vi)] $\displaystyle{0\le \psi(x,y)\le Lx^2}$ in $\mcl{Q}^f_{\eps_0}$ for some constant $L>0$;

\item[(vii)] $\displaystyle{0\le \frac{\psi_x}{x}\le \frac{2-\delta}{\gam+1}}$ in $\mcl{Q}^f_{\eps_0}$ for some constant $\delta\in(0, \frac 12)$.
\end{itemize}
Then it directly follows from Theorem \ref{theorem:1} that
\begin{itemize}
\item[(a)] $\forall y\in [0, f(0))$,
\begin{equation*}
  \lim_{x\to 0+}\psi_{xx}(x,y)=\frac{1}{\gam+1};
\end{equation*}
\item[(b)] $\displaystyle{\psi|_{\ol{\mcl{Q}^f_{\eps_0}}}}$ is not $C^2$ at the point $(0, f(0))$.
\end{itemize}

According to the statement (a), the radial derivative of the flow velocity$(=\nabla_{\rx}\Phi)$ on the pseudo-sonic arc $\Gam_0$(Fig.\ref{fig1}, Fig.\ref{fig2}(left)) is nonzero. While the velocity field $\nabla_{\rx}\Phi$ is discontinuous on the shock $\ol{S_0\cup\shock}$, it is continuous on $\Gam_0$. What the statement (a) indicates, however, is that the pseudo-sonic arc $\Gam_0$ is a {\emph{weak discontinuity}} in the sense that a derivative of the velocity field is discontinuous on $\Gam_0$. We also point out that another pseudo-sonic arc $\Gam_1$ in Fig.\ref{fig2} due to the presence of a normal shock state in $\Om_1$ is also a {\emph{a weak discontinuity}}.
Now, a question arises naturally:

{\emph{Given a second order equation with a degeneracy of Keldysh type, does the degeneracy always result in a discontinuity in a second order derivative of its solutions on the degenerate interface? }}
\bigskip

Another example introduced in the next section indicates that the answer to the question is `no'.

\section{A sonic interface as a regular interface}
\label{section:3}
Given a constant $L>0$, define
\begin{equation*}
  \Om_L:=\{{\rx}=(x_1, x_2)\in \R^2: 0<x_1<L,\,\,|x_2|<1\}.
\end{equation*}
The boundary $\der\Om_L$ consists of the entrance $\Gam_0=\{0\}\times[-1,1]$, $\Gamw:=(0,L)\times \{\pm 1\}$, and the exit $\Gamex:=\{L\}\times [-1,1]$.

For two fixed constants $\gam>1$ and $\barrhoi>0$, consider the steady Euler-Poisson system
\begin{equation}
\label{E-P}
  \begin{aligned}
  \nabla\cdot(\rho {\bf u})=0\\
  \nabla\cdot(\rho{\bf u}\otimes {\bf u})+\nabla p=\rho \nabla\Phi\\
  \nabla\cdot (\rho \mcl{B}{\bf u})=\rho{\bf u}\cdot\nabla\Phi\\
  \Delta\Phi=\rho-\barrhoi
  \end{aligned}
\end{equation}
with
\begin{equation*}
  \mcl{B}=\frac 12|{\bf u}|^2+\frac{\gam p}{(\gam-1)\rho}
\end{equation*}
for a fixed adiabatic exponent $\gam>1$.

Any one-dimensional solution $(\rho, {\bf u}, p, \Phi)=(\bar{\rho}, \bar u_1{\bf e}_1, \bar p, \bar{\Phi})(x_1)$ with $\bar{\rho}>0$ and $\bar u_1>0$ can be given as
\begin{equation*}
  \begin{aligned}
  \bar{\rho}(x_1)&=\frac{J}{\bar u_1(x_1)}\\
  \bar p(x_1)&=S_0\bar{\rho}^{\gam}(x_1)\\
  \bar{\Phi}(x_1)&=\frac{u_0^2}{2}+\frac{\gam S_0}{\gam-1}\left(\frac{J}{u_0}\right)^{\gam-1}+\int_0^{x_1}\bar E(t)\,dt
  \end{aligned}
\end{equation*}
with $(\bar u_1, \bar E)$ solving
\begin{equation}
\label{E-P:1-d}
  \begin{cases}
  \bar u_1'=\frac{\bar E\bar u_1^{\gam}}{\bar u_1^{\gam+1}-\us^{\gam+1}}\\
  \bar E'=\frac{J}{\bar u_1}-\barrhoi
  \end{cases}\quad\tx{for $x_1>0$},\quad (\bar u_1, \bar E)(0)=(u_0, E_0)
\end{equation}
for constants $S_0>0$, $J>0$, $u_0>0$ and $E_0\in \R$. Here, the constant $\us$ represents the speed at the sonic state, and is explicitly given by
\begin{equation*}
  \us=(\gam S_0 J^{\gam-1})^{\frac{1}{\gam+1}}.
\end{equation*}

Define $(\barui, \zeta_0):=(\frac{J}{\barrhoi}, \frac{\barui}{\us})$, and assume that
\begin{equation*}
  \zeta_0>1.
\end{equation*}
For $H:(0,\infty)\rightarrow \R$ given by
\begin{equation*}
  H(u):=\frac{J}{\barui}\int_{\us}^u \frac{1}{t^{\gam+1}}(t^{\gam+1}-\us^{\gam+1})(\barui-t)\,dt,
\end{equation*}
if $(\bar u_1, \bar E)(x_1)$ is a $C^1$-solution to \eqref{E-P:1-d}, then it can be directly checked from \eqref{E-P:1-d} that
\begin{equation*}
 \frac 12 \bar E^2-H(\bar u_1)=\frac 12 E_0^2-H(u_0)\quad\tx{for $x_1>0$}
\end{equation*}
as long as the solution exists.

We call a set
\begin{equation}
\label{definition:critical trajectory}
  \mcl{T}:=\{(u,E): \frac 12 E^2-H(u)=0\}
\end{equation}
{\emph{the critical trajectory}} on the $uE$-plane. Further, we call a set
\begin{equation*}
  \mcl{T}_{\rm acc}:=\{(u, E)\in \mcl{T}: (u-\us)E\ge 0\}
\end{equation*}
{\emph{the critical trajectory with an acceleration}}(Fig.\ref{fig4}).
\begin{psfrags}
\begin{figure}[htp]
\centering
\psfrag{A}[cc][][0.8][0]{$\phantom{aaaa}(u_0, E_0)$}
\psfrag{B}[cc][][0.8][0]{$u_s$}
\psfrag{C}[cc][][0.8][0]{$\barui$}
\psfrag{u}[cc][][0.8][0]{$u$}
\psfrag{E}[cc][][0.8][0]{$E$}
\includegraphics[scale=0.5]{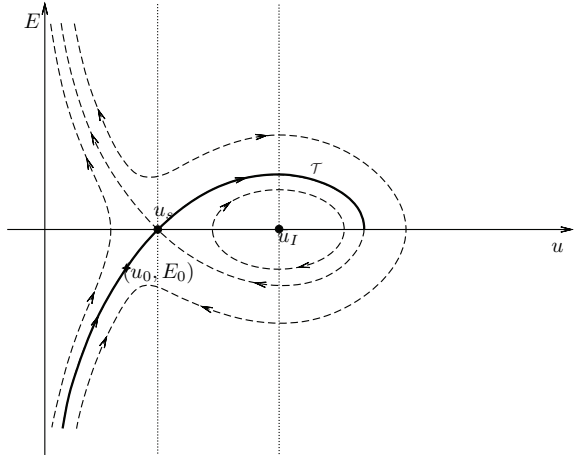}
\caption{The critical trajectory}\label{fig4}
\end{figure}
\end{psfrags}

For the initial data $(u_0, E_0)$ in \eqref{E-P:1-d}, suppose
\begin{equation}
\label{condition-initial data}
  (u_0, E_0)\in \mcl{T}_{\rm acc},\quad u_0<\us.
\end{equation}
\begin{lemma}[Lemma 1.1 in \cite{bae2023steady}]
\label{lemma:1}
The initial value problem \eqref{E-P:1-d} with $(u_0, E_0)$ satisfying \eqref{condition-initial data}  has a unique smooth solution $(\bar u_1, \bar E)$ with the following properties:
\begin{itemize}
\item[(i)] there exists a finite constant $l_{\rm max}>0$ such that
    \begin{equation*}
      \bar u_1'(x_1)>0\quad\tx{for $x_1\in[0, l_{\rm max})$};
    \end{equation*}
\item[(ii)] $\displaystyle{\lim_{x \to l_{\rm max}- }}\bar u_1'(x_1)=0$;
\item[(iii)] $\mcl{T}_{\rm acc}\cap\{(\bar u_1, \bar E)(x_1):0\le x_1\le l_{\rm max}\}=\mcl{T}_{\rm acc}\cap\{(u, E): u\ge u_0\}$ (Fig.\ref{fig4});
\item[(iv)] there exists a unique constant $\ls \in(0, l_{\rm max})$(Fig.\ref{fig5}) such that
    \begin{equation*}
      \bar u_1(x_1)\begin{cases}
      <\us\quad&\mbox{for $x_1<\ls$}\\
      =\us\quad&\mbox{at $x_1=\ls$}\\
      >\us\quad&\mbox{for $x_1>\ls$}
      \end{cases}.
    \end{equation*}
\end{itemize}
\end{lemma}

\begin{psfrags}
\begin{figure}[htp]
\centering
\psfrag{a}[cc][][0.8][0]{$\Gam_0$}
\psfrag{b}[cc][][0.8][0]{$\Gam_w$}
\psfrag{d}[cc][][0.8][0]{$\Gam_L$}
\psfrag{c}[cc][][0.8][0]{$\phantom{aaaa}x_1=l_s$}
\psfrag{e}[cc][][0.8][0]{subsonic}
\psfrag{f}[cc][][0.8][0]{supersonic}
\psfrag{g}[cc][][0.8][0]{$\phantom{aaa}\Gamma_{\rm sonic}$}
\includegraphics[scale=0.5]{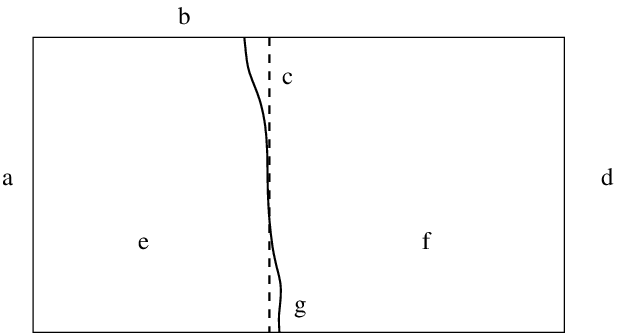}
\caption{}\label{fig5}
\end{figure}
\end{psfrags}
For the solution $(\bar u_1, \bar E)$ given in Lemma \ref{lemma:1}, we can write as
\begin{equation}
\label{E-representation}
  \bar E={\rm sgn}(\bar u_1-\us)\sqrt{2H(\bar u_1)}\quad\tx{for $x_1\in[0, l_{\rm max}]$}.
\end{equation}
Set
\begin{equation*}
  \bar{\vphi}(x_1):=\int_0^{x_1}\bar u_1(t)\,dt.
\end{equation*}
It easily follows from \eqref{E-P:1-d} and \eqref{E-representation} that
\begin{equation*}
  \left((\der_1\bar{\vphi})^{\gam+1}-\us^{\gam+1}\right)\der_{11}\bar{\vphi}-{\rm sgn}(\der_1\bar{\vphi}-\us)\sqrt{2H(\der_1\bar{\vphi})}(\der_1\bar{\vphi})^{\gam}=0\quad\tx{for $x_1\in[0, l_{\rm max}]$}.
\end{equation*}
Therefore, $\bar{\vphi}$ can be regarded as a solution to a second order equation with a degeneracy of Keldysh type, and the smoothness of $\bar{\vphi}$ indicates that the degenerate interface $x_1=l_s$ is not a weak discontinuity. So we are given with an example for a degenerate interface of a new type for a second order equation with a degeneracy of Keldysh type. Naturally, the following questions arise:
\begin{itemize}
\item[-] {\emph{Does there exist a multi-dimensional solution to \eqref{E-P} with a degenerate interface?}}
\item[-] {\emph{If so, what is the regularity of the solution across the degenerate interface?}}
\end{itemize}
Let us assume that
\begin{equation*}
  \mcl{B}-\Phi=0,
\end{equation*}
which we call {\emph{a pseudo Bernoulli's law}}, and let us write the pressure function $p$ as
\begin{equation*}
  p=S\rho^{\gam}.
\end{equation*}
For $\rho>0$ and ${\bf u}\cdot{\bf e}_1(=:u_1)>0$, the system \eqref{E-P} can be rewritten as
\begin{equation}\label{full EP rewritten}
  \left\{
\begin{aligned}
& \nabla\cdot ({\rho{\bf u}})=0 \\
&\nabla\times {\bf u}=\frac{\rho^{\gam-1}\der_{2} S}{(\gam-1)u_1}\\
&  {\rho{\bf u}}\cdot \nabla S=0\\
&\frac 12|{\bf u}|^2+\frac{\gam S \rho^{\gam-1}}{\gam-1}=\Phi\\
& \Delta \Phi=\rho-\barrhoi.
\end{aligned}%
\right.
\end{equation}

Take the one-dimensional solution $(\bar{\rho}, \bar u_1{\bf e}_1, \bar p(=S_0\bar{\rho}^{\gam}), \bar{\Phi})$ with $(\bar u_1, \bar E)$ given from Lemma \ref{lemma:1}.
For a two-dimensional velocity field ${\bf u}$, we represent as
\begin{equation*}
  {\bf u}=\nabla\vphi+\nabla^{\perp}\phi,
\end{equation*}
and define
\begin{equation*}
  (\psi,\Psi):=(\vphi, \Phi)-(\bar{\vphi}, \bar{\Phi}).
\end{equation*}
Then we can further rewrite \eqref{full EP rewritten} as a nonlinear system for $(\psi, \Psi, \phi, S)$ as follows:
\begin{align}
\label{equation for psi}
 \sum_{i,j=1}^2 a_{ij}\der_{ij}\psi+  b_1\der_1\psi+
    \bar u_1\der_1\Psi
    +(\gam-1)\bar u_1'\Psi={f}_1(S,\nabla\psi, \nabla^{\perp}\phi, D(\nabla^{\perp}\phi) \Psi, \nabla\Psi)&\\
\label{equation for Psi}
\Delta \Psi-\frac{1}{\gam S_0\bar{\rho}^{\gam-2}}\Psi+\frac{\bar u_1}{\gam S_0\bar{\rho}^{\gam-2}}\der_1\psi={f}_2(S,\nabla\psi, \nabla^{\perp}\phi, \Psi, \nabla\Psi)&\\
\label{equation for phi}
  -\Delta \phi= {f}_3(S,\nabla\psi, \nabla^{\perp}\phi, \Psi, \nabla\Psi)&\\
    \label{equation for T}
(\nabla\bar{\vphi}+\nabla\psi+\nabla^{\perp}\phi)\cdot \nabla S=0
    \end{align}
for
\begin{equation*}
  \begin{aligned}
  a_{ij}&=(\gam+1)\left(\bar{\Phi}+\Psi-\frac 12|{\bf v}|^2\right)\delta_{ij}
  -({\bf v}\cdot {\bf e}_i)({\bf v}\cdot {\bf e}_j)\,\,\tx{with ${\bf v}=\nabla\bar{\vphi}+\nabla\psi+\nabla^{\perp}\phi$}\\
  b_1&=\bar E-(\gam+1)\bar u_1'\bar u_1.
  \end{aligned}
\end{equation*}

To find a two-dimensional solution to \eqref{E-P} as a small perturbation of $(\bar{\rho}, \bar u_1{\bf e}_1, \bar p, \bar{\Phi})$. we prescribe the following boundary conditions:
\begin{equation}
\label{bc for E-P}
  \begin{aligned}
  {\bf u}\cdot {\bf e}_2=w_{\rm en},\quad S=S_{\rm en},\quad \der_1\Phi=E_{\rm en}\quad&\mbox{on $\Gam_0$}\\
  {\bf u}\cdot {\bf e}_2=0,\quad \der_2\Phi=0 \quad &\mbox{on $\Gam_w$}\\
  \Phi=\bar{\Phi}\quad&\mbox{on $\Gam_L$}
  \end{aligned}
\end{equation}
for three functions $w_{\rm en},\, S_{\rm en},\, E_{\rm en}:[-1,1]\rightarrow \R $ satisfying
\begin{equation}
\label{perturbation}
  \mcl{P}(w_{\rm en}, S_{\rm en}, E_{\rm en}):=\|w_{\rm en}\|_C^5([-1,1])+\|(S_{\rm en}, E_{\rm en})-(S_0, E_0)\|_{C^4([-1,1])}\le \sigma
\end{equation}
for some small constant $\sigma>0$, and satisfying the compatibility conditions
\begin{equation*}
\begin{aligned}
\left(\frac{d}{dx_2}\right)^kE_{\rm en}=\left(\frac{d}{dx_2}\right)^kS_{\rm en}=0\quad&\tx{at $|x_2|=1$ for $k=1,3$},\\
\left(\frac{d}{dx_2}\right)^lw_{\rm en}=0\quad&\tx{at $|x_2|=1$ for $l=0,2,4$.}
\end{aligned}
\end{equation*}
In the framework of the Helmhotz decomposition, \eqref{bc for E-P} becomes
\begin{equation}
\label{bc for decomposition}
  \begin{aligned}
  \der_2\psi=\om_{\rm en},\,\, \der_1\phi=0,\,\, S=S_{\rm en},\,\, \der_1\Psi=E_{\rm en}-E_0\quad&\mbox{on $\Gam_0$}\\
  \der_2\psi=0,\,\, \phi=0,\,\, \der_2\Psi=0\quad&\mbox{on $\Gam_w$}\\
  \phi=0,\,\, \Psi=0\quad&\mbox{on $\Gam_L$}.
  \end{aligned}
\end{equation}
The boundary condition $\phi=0$ on $\Gam_L$ is added in \eqref{bc for decomposition} for the well-posedness of the boundary value problem \eqref{equation for psi}--\eqref{equation for T} with \eqref{bc for decomposition} for $(\psi, \Psi, \phi, S)$.
As we seek for a solution $(\psi, \Psi, \phi, S)$ with $\|(\psi, \Psi, \phi)\|_{W^{1,\infty}}$ being small, we first investigate a modified equation of \eqref{equation for psi}:
\begin{equation}
\label{3-12}
\bar a_{11}\der_{11}w+\bar a_{22}\der_{22}w+b_1\der_1w=f
\end{equation}
for $(\bar a_{11}, \bar a_{22})=(a_{11},a_{22})$ with $\psi=\phi=\Psi=0$. For the normalized coefficients
\begin{equation}
\label{normalized coefficinets}
  \alp_{11}:=\frac{\bar a_{11}}{\bar a_{22}}=1-\left(\frac{\bar u_1}{\us}\right)^{\gam+1},\quad
  \beta_1:=\frac{b_1}{\bar a_{22}}=\frac{(\bar E-(\gam+1)\bar u_1'\bar u_1)\bar u_1^{\gam-1}}{\us^{\gam+1}},
\end{equation}
define a linear differential operator
\begin{equation}
\label{normalized operator}
  \mcl{L}w:=\alp_{11}\der_{11}w+\der_{22}w+\beta_1\der_1 w.
\end{equation}
According to Lemma \ref{lemma:1}(iv), the operator $\mcl{L}$ is elliptic in $\Om_L\cap\{x_1<l_s\}$, hyperbolic in $\Om_L\cap\{x_1>l_s\}$, and degenerate on $\Om_L\cap \{x_1=l_s\}$. Therefore, it is a mixed type operator with a degeneracy of Keldysh type if $L>l_s$. From this, it is easy to see that the equation \eqref{equation for psi}, as a second order equation for $\psi$, is a mixed type with a degeneracy of Keldysh type.

One of differences between the two equations \eqref{2-1} and \eqref{equation for psi} (or \eqref{3-12}) is that the equation \eqref{2-1} is hyperbolic before($x<0$) the degenerate boundary(=sonic arc), and turns to be elliptic after($x>0$) the degenerate interface, while the equation \eqref{3-12} changes its type from being elliptic to being hyperbolic across the degenerate interface. Does this difference lead to a different result on the regularity of solutions to the boundary value problem \eqref{equation for psi}--\eqref{equation for T} with \eqref{bc for decomposition}?

\begin{theorem}[Theorem 2 in \cite{bae2023steady}]
\label{theorem-HD}
Given constants $(\gam, \zeta_0, J, S_0, E_0)$ satisfying
\begin{equation*}
\gam>1,\quad \zeta_0>1,\quad S_0>0,\quad E_0<0,
\end{equation*}
suppose that $(u_0, E_0)\in \Tac$, which is equivalent to $0<u_0<\us$.
Then one can fix two constants $\bJ$ and $\ubJ$ depending only on $(\gam, \zeta_0, S_0)$ with $0<\bJ<1<\ubJ<\infty$
so that whenever the background momentum density $J(=\bar{\rho}\bar u_1)$ satisfies
\begin{center}
either $0<J\le \bJ$ or $\ubJ\le J<\infty$,
\end{center}
there exists a constant $d\in(0,1)$ depending on $(\gam, \zeta_0, S_0,J)$ so that if the two constants $u_0$ and $L$ are fixed to satisfy
\begin{equation}
\label{almost sonic condition1 full EP}
  1-d\le \frac{u_0}{\us}<1 <\frac{\bar u_1(L)}{\us}\le 1+d,
\end{equation}
and if the constant $\sigma>0$ in \eqref{perturbation} is fixed sufficiently small depending only on $(\gam, \zeta_0, S_0, E_0, J, L)$,
then the boundary value problem \eqref{equation for psi}--\eqref{equation for T} with \eqref{bc for decomposition} has a unique solution $(\psi, \phi, \Psi, S )$ that satisfies the estimate
    \begin{equation}
    \label{solution estimate HD}
    \begin{split}
    \|(\psi, \phi, \Psi, S-S_0)\|_{H^4(\Om_L)}
    \le
    C\mcl P(S_{\rm en}, E_{\rm en}, w_{\rm en})
    \end{split}
    \end{equation}
    for some constant $C>0$ depending on $(\gam, \zeta_0, S_0, E_0, J, L)$.
    \medskip

Furthermore, there exists a function $\fsonic:[-1,1]\rightarrow (0, L)$ such that
\begin{equation}
    \label{sonic boundary is a graph pt}
  \frac{|{\bf u}|}{\sqrt{\gam S\rho^{\gam-1}}}\begin{cases}
  <1\quad&\mbox{for $x_1<\fsonic(x_2)$}\\
  =1\quad&\mbox{for $x_1=\fsonic(x_2)$}\\
  >1\quad&\mbox{for $x_1>\fsonic(x_2)$}
  \end{cases},
    %\begin{split}
    %|{\bf u}|<\sqrt{\gam S\rho^{\gam-1}}\quad&\mbox{for $x_1<\fsonic(x_2)$\quad(subsonic)},\\
    %|{\bf u}|=\sqrt{\gam S\rho^{\gam-1}}\quad&\mbox{for $x_1=\fsonic(x_2)$\quad(sonic)},\\
    %|{\bf u}|>\sqrt{\gam S\rho^{\gam-1}}\quad&\mbox{for $x_1>\fsonic(x_2)$\quad(supersonic)}.
     % \end{split}
    \end{equation}
and the function $\fsonic$ satisfies
\begin{equation}
    \label{estimate of sonic boundary pt}
      \|\fsonic-\ls\|_{H^2((-1,1))}+\|\fsonic-\ls\|_{C^1([-1,1])}\le C\mcl P(S_{\rm en}, E_{\rm en}, w_{\rm en})
    \end{equation}
    for some constant $C>0$ depending on $(\gam, \zeta_0, S_0, E_0, J, L)$.
\end{theorem}

We are given from Theorem \ref{theorem-HD} with a classical solution $(\rho, {\bf u}, p)$ of \eqref{E-P} with the sonic interface $x_1=\fsonic(x_2)$ across which the velocity field ${\bf u}$ is not only continuous, but also its derivative is continuous. In other words, $x_1=\fsonic(x_2)$ is a degenerate interface but not a weak discontinuity. We call this sonic interface as a {\emph{regular interface}} in the sense that the solution $(\rho, {\bf u}, p)$ is a classical solution across the sonic interface.

\section{Discussion}
\label{section:discussion}
In \S\ref{section:2}, we show an example of a (pseudo) sonic interface as a weak discontinuity in the sense that a velocity field is continuous but its derivative is discontinuous on the interface. Such an example is given from a solution to a mixed-type equation with a degeneracy of Keldysh type. More precisely, the type of the equation changes from being hyperbolic to being elliptic through a degenerate interface of codimension one. In \S\ref{section:3}, we show an example of a sonic interface as a regular interface on which a velocity field is $C^1$. This example is also given from a solution to a degenerate equation of Keldysh type, where the equation changes its type from being elliptic to being hyperbolic.
\medskip

The interesting point is that the equations given in \S\ref{section:2}--\ref{section:3} are both Keldysh type, and that the lower order derivative terms `$(-1+O_4)\psi_x$' from \eqref{2-1}, and `$b_1\der_1\psi$' from \eqref{equation for psi} have significant contributions to Theorem \ref{theorem:1} and Theorem \ref{theorem-HD}, respectively. In spite of those similarities, the equations \eqref{2-1} and \eqref{equation for psi} give different regularity results. Why?

\subsection{The existence of a special smooth solution}
In \S\ref{section:3}, the classical solution of \eqref{E-P} with a sonic interface as a regular interface(Theorem \ref{theorem-HD}) is constructed as a small perturbation of a smooth one-dimensional solution(Lemma \ref{lemma:1}). Therefore, one may suggest that the existence of a special smooth solution can be a clue to determine a destiny of a sonic interface in general.
\medskip

It is true that a smooth one-dimensional solution comes in handy to reformulate \eqref{E-P} into \eqref{equation for psi}--\eqref{equation for T} so that an iteration method can be applied. But still, it is another matter to establish the well-posedness of the boundary value problem \eqref{equation for psi}--\eqref{equation for T} with \eqref{bc for decomposition}. According to the work in \cite{bae2023steady}, the essential ingredient used to prove Theorem \ref{theorem-HD} is the strictly increasing property of $\bar u_1$ stated in Lemma \ref{lemma:1}(i). This monotonicity property combined with the method developed in \cite{KZ} yields the well-posedness of a linearized boundary value problem derived from \eqref{equation for psi}--\eqref{equation for T} with \eqref{bc for E-P}, and this yields Theorem \ref{theorem-HD}.

\medskip

\begin{figure}[htp]
\centering
\begin{psfrags}
\psfrag{a}[cc][][0.7][0]{}
\psfrag{f}[cc][][0.7][0]{$S_0$}
\psfrag{b}[cc][][0.7][0]{$\Om_0$}
\psfrag{e}[cc][][0.7][0]{$\Gam_0$}
\psfrag{c}[cc][][0.7][0]{$\Om$}
\psfrag{d}[cc][][0.7][0]{$\phantom{a}\theta_w=\frac{\pi}{2}$}
\includegraphics[scale=0.8]{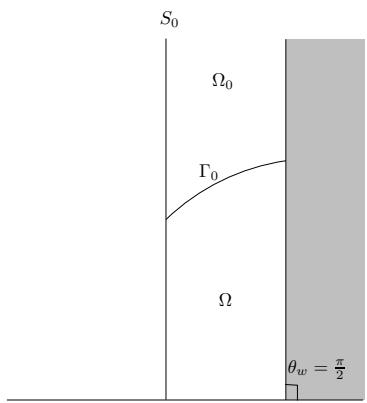}
\caption{Self-similar normal shock reflection ($S_0$: reflected normal shock, $\vphi=\vphi_0$ in $\Om_0\cup\Gam_0\cup\Om$, see \cite{chen2010global})}
\label{fig6}
\end{psfrags}
\end{figure}
But, there is another example that shows that the existence of a special smooth solution does not necessarily tells the destiny of a sonic interface in general. As repeatedly pointed out, the sonic arc $\Gam_0$ in a self-similar regular shock reflection is a weak discontinuity for $\theta_w\in (\theta_{\rm sonic}, \frac{\pi}{2})$(see \S\ref{section:2}). For the wedge-angle $\theta_w=\frac{\pi}{2}$, however, the sonic arc is not a weak discontinuity(Fig. \ref{fig6}) because the state behind the reflected normal shock $S_0$ is simply given by $\vphi=\vphi_0$ for a quadratic polynomial function $\vphi_0$ of self-similar variables $\bmxi=(\xi_1, \xi_2)$ in the form of \eqref{pseudo potential0}. In \cite{chen2010global}, a self-similar regular shock reflection solution for $\theta_w\in(\frac{\pi}{2}-\sigma, \frac{\pi}{2})$ is constructed as a small perturbation of the normal shock reflection solution, but the nature of the sonic arc $\Gam_0$ for $\theta_w<\frac{\pi}{2}$ is not determined by the arc $\Gam_0$ for $\theta_w=\frac{\pi}{2}$.

\subsection{Acceleration VS. deceleration in a transonic transition}
Notice that the self-similar flow(Fig.\ref{fig1} and Fig.\ref{fig2}(left)) introduced in \S\ref{section:2} turns from being (pseudo) supersonic in $\Om_0$ to being (pseudo) subsonic in $\Om$ through the (pseudo) sonic arc $\Gam_0$. This is a decelerating transition on $\Gam_0$, and Theorem \ref{theorem:1} implies that $\Gam_0$ is a weak discontinuity. On the other hand, the $C^1$-transonic solution of Euler-Poisson system \eqref{E-P} given by Theorem \ref{theorem-HD} in \S\ref{section:3} has an accelerating transition on the sonic interface $x_1=\fsonic(x_2)$. \begin{question}
Does a decelerating transonic transition always yield a sonic interface as a weak discontinuity, while accelerating transonic transition does not?

\end{question}

Here is another interesting example to investigate. For the set $\mcl{T}$ given by \eqref{definition:critical trajectory}, define
\begin{equation*}
  \mcl{T}_{\rm dec}:=\{(u, E)\in \mcl{T}: (u-\us)E\le 0\}.
\end{equation*}
We call the set $\mcl{T}_{\rm dec}$ {\emph{the critical trajectory with a deceleration}}. By modifying the proof of Lemma \ref{lemma:1} given in \cite{bae2023steady}, one obtains the following lemma:
\begin{lemma}
\label{lemma:2}
For the initial data $(u_0, E_0)$ in \eqref{E-P:1-d}, suppose
\begin{equation*}
  (u_0, E_0)\in \mcl{T}_{\rm dec},\quad u_0>\us.
\end{equation*}
Then, the initial value problem \eqref{E-P:1-d} has a unique smooth solution $(\bar u_1, \bar E)$ with the following properties: \\
(Case 1) For $1<\gam<2$,
\begin{itemize}
\item[(i)] there exists a finite constant $\til l_{\rm max}>0$ such that
    \begin{equation*}
      \bar u_1'(x_1)<0\quad\tx{for $x_1\in[0, \til l_{\rm max})$};
    \end{equation*}
\item[(ii)] $\displaystyle{\lim_{x_1 \to \til l_{\rm max}- }\bar u_1'(x_1)=0\,\,\tx{and}\,\,\lim_{x_1 \to \til l_{\rm max}- }\bar E(x_1)=\infty}$;
\item[(iii)] $\mcl{T}_{\rm dec}\cap\{(\bar u_1, \bar E)(x_1):0\le x_1\le \til l_{\rm max}\}=\mcl{T}_{\rm dec}\cap\{(u, E): u\le u_0\}$ (Fig.\ref{fig4});
\item[(iv)] there exists a unique constant $\til l_s \in(0, \til l_{\rm max})$ such that
    \begin{equation*}
      \bar u_1(x_1)\begin{cases}
      >\us\quad&\mbox{for $x_1<\til l_s $}\\
      =\us\quad&\mbox{at $x_1=\til l_s $}\\
      <\us\quad&\mbox{for $x_1>\til l_s $}
      \end{cases}.
    \end{equation*}
\end{itemize}
(Case 2) For $\gam\ge 2$, the statements given in (i)--(iv) hold with $\til l_{\rm max}=\infty$.
\end{lemma}

The smooth one-dimensional solution of the Euler-Poisson system \eqref{E-P} has a decelerating speed and a sonic interface at $x_1=\til l_s$. Similarly to Theorem \ref{theorem-HD}, would it be possible to establish the existence of a classical solution to \eqref{E-P} as a small perturbation of this one-dimensional decelerating smooth transonic solution? In proving Theorem \ref{theorem-HD}, the key property used in \cite{bae2023steady} is that there exists a constant $\lambda_L>0$ satisfying
\begin{equation}
\label{KZ inequality}
  -2\beta_1-(2m-1)\der_1\alp_{11}\ge \lambda_L\quad\tx{in $\Om_L$ for $m=0,1,2,3$}
\end{equation}
for the coefficients $\alp_{11}$ and $\beta_1$ given by \eqref{normalized coefficinets}. This inequality enables to apply \cite[Theorem 1.7]{KZ} to achieve a priori $H^{m+1}$ estimates of solutions obtained from Theorem \ref{theorem-HD} for $m=0,1,2,3$. More precisely, a direct computation yields the representation
\begin{equation*}
  -2\beta_1-(2m-1)\der_1\alp_{11}= \frac{\bar u_1'}{\us^{\gam+1}}\left(2m(\gam+1)\bar u_1^{\gam}+(\gam-1)\bar u_1^{\gam}+2\frac{\us^{\gam+1}}{\bar u_1}\right).
\end{equation*}
And, the inequality \eqref{KZ inequality} is obtained by Lemma \ref{lemma:1}.
\medskip

For one-dimensional smooth transonic solution $(\bar u_1, \bar E)$ lying on the critical trajectory $\mcl{T}_{\rm dec}$ with a deceleration, Lemma \ref{lemma:2} implies that the inequality \eqref{KZ inequality} does not hold because $\bar u_1'<0$. So one cannot apply \cite[Theorem 1.7]{KZ} to establish the existence of a classical solution to \eqref{E-P} as a small perturbation of $(\bar u_1, \bar E)$. This opens to a possibility that a multi-dimensional solution to \eqref{E-P} given as a small perturbation of $(\bar u_1, \bar E)$ with $\bar u_1'<0$ may contain a sonic interface as a weak discontinuity. This is an open problem to be investigated in the future.

\medskip

\vspace{.25in}
\noindent
{\bf Acknowledgements:} The author shows her sincere respect to Prof. Gui-Qiang Chen's for his ambition and contribution to mathematics.
The research of Myoungjean Bae was supported in part by  Samsung Science and Technology Foundation under Project Number SSTF-BA1502-51.

\bigskip
\bibliographystyle{acm}
\bibliography{References_EP_smooth_tr}

\end{document}